\documentclass[12pt]{article}
\usepackage{amsfonts}
\usepackage{amssymb}
\usepackage{amsmath}
\usepackage{amsxtra}
\usepackage{amscd}
\usepackage{amsthm}
\usepackage{eucal}

\title{Interpreting Random Hypergraphs in Pseudofinite Fields}
\author{\"{O}zlem Beyarslan \\ University of Illinois at Chicago}
\date{\today}
\newtheorem{theorem}{Theorem}

\newtheorem{proposition}[theorem]{Proposition}
\newtheorem{lemma}[theorem]{Lemma}
\newtheorem{corollary}[theorem]{Corollary}
\newtheorem{fact}[theorem]{Fact}
\newtheorem{notation}[theorem]{Notation}

\def\nn {\mathbb{N}}

\def\qq {\mathbb{Q}}

\newfont{\mym}{msbm10}

\def\1{^{\circ}}
\def\2{_{\circ}}

\def\claim{\noindent {\bf Claim:} $\,$}
\def\proof{\noindent {\bf Proof:} $\,$}

\def\-#1{\overline{#1}}
\def\a{\alpha}
\def\b{\beta}

\def\ontor{\rlap{$\to$}\kern.05pt\to}

\def\qed{$\Box$}

\def\ch{\mathop{\rm char}\nolimits}
\def\deg{\mathop{\rm deg}\nolimits}

\def\Sym{\mathop{\rm Sym}\nolimits}
\def\Alt{\mathop{\rm Alt}\nolimits}
\def\Fix{\mathop{\rm Fix}\nolimits}
\def\Aut{\mathop{\rm Aut}\nolimits}

\def\lcm{\mathop{\rm lcm}\nolimits}

\def\Gal{\mathop{\rm Gal}\nolimits}

\def\Id{\mathop{\rm Id}\nolimits}
\def\mod{\mathop{\rm mod}\nolimits}

\def\Id{\mathop{\rm Id}\nolimits}

\def\mod{\mathop{\rm mod}\nolimits}
\def\qed{\hfill$\Box$}
\def\proof{\noindent {\bf Proof:} $\,$}
\def\notation{\noindent {\bf Notation:} $\,$}

\begin{document}
\maketitle

\section{Introduction}
This article gives a positive answer to a question posed by D.\
Macpherson (Ravello 2002 \cite{qm}, Question 14):

``Can we interpret the random $n$-ary hypergraph in a pseudofinite
field?"

A {\em pseudofinite field} is an infinite field that satisfies all
first-order sentences that hold in {\bf every} finite field. An
example of a pseudofinite field is an infinite ultraproduct of
finite fields. The theory of pseudofinite fields was first studied
by J.\ Ax  in his 1968  article ``The elementary theory of finite
fields". In this article, among other results, Ax proves that a
field $F$ is pseudofinite if and only if it is perfect, has a
unique extension of degree $n$ for every $n\in \nn^{>0}$ and is
pseudo algebraically closed (PAC), that is, every absolutely
irreducible variety defined over $F$ has an $F$-rational point.

In 1980 J.\ L.\ Duret showed \cite {Duret} that the theory of
pseudofinite fields is unstable, as the random graph is definable:
given a pseudofinite field $F$ of characteristic different from $2$,
put an edge between any two distinct points in $F$ in case their sum
is a square in $F$.

During the early 1990's Hrushovski \cite{pac} showed that the theory
of pseudofinite fields, although unstable, is not so ``bad" in the
sense that, some of the methods from stability theory can still be
applied here.

An {\em $n$-hypergraph} is a graph whose edges, instead of
connecting just two vertices, connect $n$ distinct vertices. A {\em
random $n$-hypergraph} on a set $A$ is a tuple $(A,\, H)$ where $H$
is a subset of $A^{[n]}$ satisfying the following sentence for every
$m$ and $k$: for all $a_1,\ldots,a_m$ and $b_1,\ldots,b_k$ in
$A^{[n-1]}$, distinct, there is an element $c\in A$, such that
$a_1\cup \{c\},\ldots,a_m\cup \{c\} \in H$  and $b_1\cup
\{c\},\ldots,b_k\cup \{c\} \not\in H$.

Hrushovski proved in \cite{pac} that it is not possible to interpret
a random $(n+1)$-ary hypergraph in a random $n$-ary hypergraph. This
proves that the complexity of the random $n$-ary graphs strictly
increases with $n$.

\vskip .1in {\bf Acknowledgements.} Many thanks to Zo\'e Chatzidakis
and Ehud Hrushovski for their extremely valuable contributions. They
suggested some of the approaches and ideas on which this paper is
built. I am truly grateful to both of them for their hospitality and
support during my stays in Paris and Jerusalem. Also thanks to my
advisor David Marker for his continuous support during my years as a
graduate student at UIC. Finally, many thanks to Matthias
Aschenbrenner and Rahim N. Moosa for their comments on a preliminary
version of this paper.

\section{Preliminaries}

Throughout the article all the fields we will consider will be contained in a fixed algebraically closed field $\Omega$ and  $K$ will stand for
a field contained in $\Omega$.

Let $\widetilde{K}$ denote the algebraic closure of $K$ in
$\Omega$ and we will denote $\Aut(\widetilde{K}/K)$ as $\Gal(K)$.
We call a field extension $L$ of $K$ a {\em regular extension} if
$K$ is algebraically closed in $L$, i.e.\ if $\widetilde{K} \cap L
= K$.

By {\em a valuation} we will mean a real discrete valuation. A valuation of a function field $K(x)$ whose maximal ideal is generated by $f(x)$
will be denoted by $v_f$. In addition to the valuations on $K(x)$ given by the maximal ideals of $K[x]$, there is one more valuation
$v_{\infty}$ of $K(x)$ which is defined by: $v_\infty(f/g) =\deg(g)-\deg(f)$ for $f$ and $g$ in $K[t]$ where $\deg(f)$ denotes the degree of the
polynomial $f$. We will denote valuations by the letters $v$ and $w$.

Let $K'$ be a finite algebraic field extension of $K$ and $\hat{v}$ be an extension of $v$ to $K'$, i.e.\ $\hat{v}$ is a valuation of $K'$ whose
valuation ring intersected with $K$ gives the valuation ring of $v$. By $r(\hat{v}:v)$ we denote the {\em ramification index} of $\hat{v}$ over
$v$ i.e.\ $r(\hat{v}:v)$ is the unique positive integer such that for all $a\in K$ we have $\hat{v}(a)=r(\hat{v}:v)v(a)$. The residue degree of
$\hat{v}$ over $v$ is the field degree of the residue field of $\hat{v}$ over the residue field of $v$ and it is denoted by $d(\hat{v}:v)$. Note
that if the residue field of $v$ is algebraically closed $d(\hat{v}:v)=1$ for every extension $\hat{v}$ of $v$.



Let $\hat{v}$ be one of the (finitely many) valuations on $K'$ that
extend $v$. Then $\hat{v}$ is said to be {\em ramified} over $v$ (or
over $K$) if $r(\hat{v}:v)>1$ and $v$ is {\em ramified} in $K'$ if
it has at least one ramified extension $\hat{v}$ to $K'$.

For a polynomial $f(X) \in K[X]$ and one of its roots $x\in
\Omega$, we will call the field extension $K(x)$ of $K$ a {\em
root field} of $f(X)$.


We call an element $\sigma$ of the absolute Galois group
$\Gal(K)$, a {\em topological generator} of $\Gal(K)$ if $\sigma$
satisfies one of the following equivalent conditions: (i) For any
finite Galois extension $L$ of $K$, $\sigma|_L$ generates
$\Gal(L/K)$. (ii) The subfield of $\widetilde{K}$ fixed by
$\sigma$ is $K$. (iii) $\langle \sigma \rangle$ is dense in
$\Gal(K)$.

It is easy to prove that $K$ has a unique extension of degree $n$
for every $n$ if and only if $\Gal(K)\simeq \widehat{\mathbb{Z}}$,
the profinite completion of $\mathbb{Z}$, and hence $\Gal(K)$ has
a topological generator. In particular the absolute Galois group
of a pseudofinite field is the profinite cyclic group
$\widehat{\mathbb{Z}}$.

\begin{proposition} \label{M} Suppose $K$ is a perfect field with exactly one extension of
degree $n$ for every positive integer $n$. Let $\sigma$ be a
topological generator of $\Gal(K)\simeq \widehat{\mathbb{Z}}$.
Suppose $L$ is a regular extension of $K$. Let $\tau \in \Gal(L)$ be
an automorphism of $\widetilde{L}$ extending $\sigma$. Let $M$ be
the subfield of $\widetilde{L}$ fixed by $\tau$. Then $\tau$ is a
topological generator of $\Gal(M)\simeq \widehat{\mathbb{Z}}$ and
$K$ is algebraically closed in $M$.
\end{proposition}

\proof Since $\tau$ extends $\sigma$, $K$ is algebraically closed in
$M$. From this it follows that if $K_n$ is the unique field
extension of $K$ of degree $n$ then $MK_n$, the join of the fields
$M$ and $K_n$ is the unique field extension of $M$ of degree $n$.
This proves the proposition.\qed

\

The next theorem characterizes the concept of elementary equivalence
of pseudofinite fields. (For model theoretical concepts we refer to
\cite{marker}).

\begin{theorem}[\cite{Ax}] \label{el} Let $E$ and $F$ be two
pseudofinite fields containing a common subfield $K$. Then $E$ and
$F$ are elementarily equivalent over $K$ if and only if the
algebraic closures of $K$ in $E$ and $F$ are isomorphic over $K$.
In particular, if $E$ is algebraically closed in $F$ then $E$ is
an elementary extension of $F$. \end{theorem}

The following proposition follows from Lemma 20.2.2 of \cite{FJ}.

\begin{theorem}[\cite{FJ}]  \label{elementary}
Let $E$ be a perfect field with at most one extension of degree
$n$ for every $n$. Then there exists a pseudofinite field $F$
containing $E$ in which $E$ is algebraically closed.
\end{theorem}

\subsection{Linearly Disjoint Extensions}

Let $E$ and $F$ be two field extensions of $K$. The fields $E$ and
$F$ are said to be {\em linearly disjoint over} $K$ if any $e_1,
\ldots, e_n \in E$ which are linearly independent over $K$ are also
linearly independent over $F$. Although not obvious from the
definition, this concept is symmetric  in $E$ and $F$ \cite[Lemma
2.5.1]{FJ}.

\begin{fact}[Lemma 2.5.2 of \cite{FJ}] \label{lin dis} Let $E$ and $F$ be two field extensions of
$K$ with $F/K$ Galois. Then $E$ and $F$ are linearly disjoint over
$K$ if and only if $E\cap F \neq K$. \end{fact}

\begin{corollary}\label{simple} Let $E$ and $F$ be two nontrivial finite
extensions of $K$ with $F/K$ Galois and $\Gal(F/K)$ simple. Then
$F$ and $E$ are linearly disjoint over $K$ if and only if $F$ is
not contained in the Galois closure of $E$ over $K$.
\end{corollary}

\proof Suppose $E$ and $F$ are not linearly disjoint over $K$.
Then by Fact \ref{lin dis}, $E\cap F > K$. Let $F_1$ be the Galois
closure over $K$ of $E\cap F$. Clearly $F_1$ is contained in the
Galois closure of $E$ over $K$. Also, since $F$ is Galois over
$K$, $F_1\leq F$. Then $\Gal(F/F_1)$ is a proper normal subgroup
of $G$. Since $G$ is simple, $\Gal(F/F_1)=\{\Id\}$ and so $F=F_1$.
That is, $F=F_1$ is contained in the Galois closure of $E$ over
$K$.

The other direction is clear from Fact~\ref{lin dis}. \qed


\vskip 0.1in

Field extensions $E_1,\ldots,E_n$  of $K$ are said to be {\em
linearly disjoint over} $K$ if each $E_i$ is linearly disjoint over
$K$ from the join of the others, equivalently if $E_i$ is linearly
disjoint from $E_1\cdots E_{i-1}$ over $K$ for every $i=2,\ldots,n$.

\begin{fact}[Lemma 2.5.6 of \cite{FJ}] \label{lindis} Let $L_1,\ldots,L_n$ be a linearly disjoint
family of Galois extensions
of $K$. Then $\Gal(L_1\ldots L_n/K) \simeq
\overset{n}{\underset{i=1}{\prod}} \Gal(L_i/K)$.
\end{fact}


\begin{lemma}\label{setwise} Finitely many distinct Galois extensions of $K$
whose Galois groups over $K$ are nonabelian finite simple groups
are linearly disjoint over $K$.
\end{lemma}

\proof Let $E_i$ ($i=1,\ldots,n$) be the Galois extensions as in the
statement of the Lemma and let $\Gal(E_i/K)=S_i$. It is enough to
show that $E_i$ is linearly disjoint from $E_1 \ldots E_{i-1}$ for
$1<i\leq n$. The claim for $i=2$ follows from Fact \ref{lin dis}.
Assuming that the claim holds for $i=n-1$, we will show that it
holds for $i=n$.

$\Gal(E_1 \cdots E_{n-1}/K) = S_1\times\ldots\times S_{n-1}$ by
induction hypothesis and Fact \ref{lindis}. Suppose for a
contradiction that $E_1 \cdots E_{n-1}$ and $E_n$ are not linearly
disjoint over $K$. Since $\Gal(E_n/K)$ is simple, by
Corollary~\ref{simple}, $E_n$ is contained in $E_1 \cdots
E_{n-1}$.

Now consider $\Gal(E_1\cdots E_{n-1}/E_n)$, which is a normal
subgroup of $S_1\times\ldots\times S_{n-1}$. By an elementary
lemma on the product of simple groups, $$\Gal(E_1\cdots
E_{n-1}/E_n) = T_1\times\ldots\times T_{n-1},$$ where $T_i$ is
either $\{1\}$ or $S_i$ for all $i=1,\ldots n-1$.  Then
$$S_1/T_1 \times \ldots \times S_{n-1}/T_{n-1} \simeq \Gal(E_1\cdots E_{n-1}/K)/
\Gal(E_1 \cdots E_{n-1}/E_n) \simeq \Gal(E_n/K)= S_n.$$ Simplicity
of $S_n$ implies that $S_k\simeq S_n$, $T_k = 1$ for some
$k=1,\ldots,n-1$ and that $S_i=T_i$ for all $i\neq k,\, n$. Thus
$$\begin{array}{lll} \Gal(E_1\cdots E_{n-1}/E_n)&=& T_1\times\ldots\times T_{n-1} = S_1\times\ldots \times
S_{k-1}\times \{1\} \times S_{k+1}\times\ldots\times S_{n-1}\\ & = &
\Gal(E_1\cdots E_{n-1}/E_k). \end{array}$$ and therefore, by the
fundamental theorem of Galois theory, $E_n=E_k$, contradicting the
assumption. \qed

\vskip .1in The lemma above still holds (with the same proof) if one
of the extensions is still simple but abelian. Observe that the only
abelian quotients of a non-abelian simple group are trivial.

\subsection{Regular Extensions}
\notation Let $f(X)$ be a polynomial in $K[X]$. For a field $F$
containing $K$, let $L$ be the splitting field of the polynomial
$f(X)$ over $F$, the Galois group $\Gal(L/F)$ is sometimes denoted
as $\Gal(f(X),F)$.

\begin{lemma}[2.6.11 of \cite{FJ}] \label{acts} Let $f(X,T)\in K[X,T]$ be a polynomial over $K$ and
$\Gal(f(X,T),K(X))$ be its Galois group over $K(X)$. Then the
polynomial $f(X,T)$ is absolutely irreducible over $K$ if and only
if $L$ is a regular extension of $K$, and in this case
$\Gal(L/K(X))$ acts transitively on the roots of $f(X,T)$ over
$K(X)$.
\end{lemma}

\begin{lemma}\label{transitive} Let $G$ be a finite group acting transitively on a finite set $A$, $|A|\geq 2$.
Then there is an element $g\in G$ such that $g(x)\neq x$ for every
$x\in A$.
\end{lemma}

\proof For $x\in A$ let $G_x$ be the stabilizer of $x$. Since $G$
acts transitively on $A$,  $[G:G_x]=|A|$ and all stabilizers are
conjugate. Any two subgroups of $G$ contain at least the identity in
their intersection, hence the cardinality of $\bigsqcup_{x\in A}G_x$
is less than $|G|$. Any $g \in G\backslash \bigsqcup_{x\in A}G_x$
will satisfy the desired condition.\qed

\vskip .1in By Lemma \ref{transitive} and Fact \ref{acts} we obtain
the following corollary.

\begin{corollary} Let $f(X,T)$ be absolutely irreducible over $K$. Let $L$ be the splitting field of $f(X,T)$ over $K(X)$.
Then there is an element $\mu$ in the Galois group $\Gal(L/K(X))$ that moves all the roots of $f(X,T)$.
\end{corollary}

\begin{fact}[2.3.11 of \cite{FJ}] \label{regular} Let $g(T)\in K[T]$ be a polynomial and $X$ an
indeterminate. Then $h(X,T)=g(T)-X\in K(X)[T]$ is absolutely
irreducible over $K$. Therefore a root field of $h(X,T)$ over
$K(X)$ is a regular extension of $K$. \end{fact}

\begin{fact}[\cite{FJ}] \label{stable} Let $f(X,T)$ be a polynomial in $K[X,T]$, separable in $T$.
Suppose $\Gal(f(X,T),K(X))\simeq \Gal(f(X,T),K_s(X))$ where $K_s$
is the separable closure of $K$. Then for any field extension $F$
of $K$ we have $\Gal(f(X,T),K(X))\simeq \Gal(f(X,T),F(X))$ and the
splitting field $L$ of $f(X,T)$ over $F(X)$ is regular over $F$.
\end{fact}

\begin{lemma}[Ch.III Corr.6, p.58 \cite{lang2}]\label{reg} Let $X\in \Omega$ be transcendental over $K$.
Let $L_1$ and $L_2$ be two algebraic extensions of $K(X)$ which
are linearly disjoint over $K(X)$ and which are regular extensions
of $K$. Then $L_1L_2$ is also a regular extension of $K$.
\end{lemma}


\vskip .2in
\subsection{Random Graphs and Hypergraphs} The theory of the
random graph is axiomatized by the statements that express the following for all natural numbers $n$ and $m$: ``for all distinct $(n+m)$
elements $x_1,\,\ldots,x_n,\, y_1,\,\ldots, y_m$ there is a $z$ such that $R(z,x_i)$ for $i=1,\ldots ,n$ and $\neg R(z,y_j)$ for $j=1,\ldots,m$.
This theory is $\omega$-categorical and has quantifier elimination.

For any set $X$, let $X^{[n]}$ denote the set of subsets of $X$ whose elements have precisely $n$ members. Then an $n$-{\em hypergraph over} $X$
is a tuple $(X,R)$ where $R$ is a subset of $X^{[n]}$. A $n$-hypergraph $(X, R)$ is called {\em random} if for every distinct $a^1,\ldots,a^m
\in X^{[n-1]}$ and for every subset $I$ of $\{1,...,m\}$ there is an element $c\in X$ such that $a^i\cup\{c\} \in R$ if and only if $i\in I$.

The countable random n-hypergraph can be constructed as the Fraiss\'e limit of  finite n-hypergraphs, hence its fist order theory is
$\omega$-categorical and has quantifier elimination by \cite[Thm 7.4.1 ]{Ho}.

Note that one can define a random $m$-hypergraph in a random
$n$-hypergraph by setting the first $n-m$ entries of the random
$n$-hypergraph to be equal to a constant $c$. On the contrary, it
was proved in \cite{pac} that if $(\omega, R)$ is isomorphic to the
random $n$-hypergraph, then $R$ is not a finite Boolean combination
of $(n-1)$-ary relations. This, together with the elimination of
quantifiers, implies that we cannot interpret a random
$n$-hypergraph in a random $(n-1)$-hypergraph.

\subsection{Symmetric Polynomials}\label{symmetric} Throughout this subsection
we let $A$ denote a commutative ring with identity and
$t_1,\ldots,t_n$ algebraically independent elements over $A$. Let
$s_{n,1},\ldots,s_{n,n}$ be the elementary symmetric polynomials in
$t_1,\ldots,t_n$ of degree $1,\ldots,n$ respectively. Thus
$$\prod_{i=1}^n (X-t_i)=X^n-s_{n,1} X^{n-1}+s_{n,2} X^{n-2}-\ldots +(-1)^ns_{n,n}. \leqno{(1)}$$
It is well known that $s_{n,1},\ldots, s_{n,n}$ form an
algebraically independent basis for the ring of symmetric
polynomials in $A[t_1,\ldots,t_n]$.



We now define $S_{n, i}$ to be the sum of all monomials of degree
$i$ over the variables $t_1,\ldots, t_n$:
$$S_{n, i}(t_1,\ldots,t_n)=\sum_{r_1+\ldots+r_n=i} t_1^{r_1}\ldots
t_n^{r_n}.$$ The polynomials $S_{n, i}$ are called the {\em
complete symmetric polynomials} in $t_1,\ldots, t_n$. The next
fact is from \cite[Section 6.1]{fulton}.

\begin{fact} \label{gen} For every
$k\leq n$, $t=(t_1,\ldots,t_n)$,
$$S_{n, k}(t) - s_{n, 1}(t)S_{n,k-1}(t)+ s_{n, 2}(t)S_{n,k-2}(t)-\ldots +(-1)^k s_{n, k}(t)=0.$$
\end{fact}

\begin{lemma}\label{s} The polynomials $S_{n, 1},\ldots, S_{n, n}$ form a
basis for the ring of symmetric polynomials in $A[t_1,\ldots,
t_n]$, that is $A[s_{n, 1},\ldots,s_{n, n}] =A[S_{n,
1},\ldots,S_{n, n}]$. \end{lemma}

\proof Obviously $A[S_{n, 1},\ldots,S_{n, n}]\leq A[s_{n,
1},\ldots,s_{n, n}]$. To prove the converse we will show that for
every $k<n$, $A[s_{n,1},\ldots,s_{n,k}] \leq
A[S_{n,1},\ldots,S_{n,k}]$ by induction on $k$. For $k=1$ there is
nothing to prove. Assume $A[s_{n,1},\ldots,s_{n,k-1}] \leq
A[S_{n,1},\ldots,S_{n,k-1}]$. It is enough to show that $s_{n,k} \in
A[S_{n,1},\ldots,S_{n,k}]$. By Fact \ref{gen} we see that
$s_{n,k}(t)$ can be written in terms of $S_{n,1}(t), \ldots,
S_{n,k}(t)$ and $s_{n,1}(t),\ldots, s_{n,k-1}(t)$. The desired
result follows by the induction hypothesis. $\hfill{\Box}$

\vskip .1 in

\notation Since the polynomial
$S_{n,n-1}(t_1,\ldots,t_n)=\sum_{r_1+\ldots+r_n=n-1} t_1^{r_1}\ldots
t_n^{r_n}$ will be used several times, we will shorten it as $S$.

\begin{lemma}\label{sim} If $a=\{a_1,\ldots,a_{n-1}\},\, b=\{b_1,\ldots,b_{n-1}\}$ are in
$F^{[n-1]}$, then
$$S(a_1,\ldots, a_{n-1},X)= S(b_1,\ldots, b_{n-1},X)$$
if and only if $a=b$. \end{lemma}

\proof Note that $S(a_1,\ldots, a_{n-1},X)= \sum_{i=0}^{n-1}
S_{n-1,i}(a_1,\ldots, a_{n-1})X^{n-1-i}$. Therefore $S(a_1,\ldots,
a_{n-1},X) = S(b_1,\ldots, b_{n-1},X)$ if and only if
$S_{n-1,i}(a_1,\ldots, a_{n-1}) = S_{n-1,i}(b_1,\ldots, b_{n-1})$
for all $i\leq n-1$. By Lemma \ref{s}, each of
$s_{n-1,1},\ldots,s_{n-1,n-1}$ can be expressed uniquely in terms of
the basis $S_{n-1,1},\ldots,S_{n-1,n-1}$. This implies that
$s_{n-1,i}(a_1,\ldots, a_{n-1}) = s_{n-1,i}(b_1,\ldots, b_{n-1})$
for all $i\leq n-1$. Hence by the fundamental equality for the
symmetric functions given in (1) above, we conclude that
$\{a_1,\ldots, a_{n-1}\}=\{b_1,\ldots, b_{n-1}\}$. $\hfill{\Box}$
\vskip 0.1in

\notation Let $f(X_1,\ldots,X_n)$ be a symmetric polynomial in
$k[X_1,\ldots,X_n]$. Let $a=\{a_1,\ldots,a_n\}$ be in $K^{[n]}$.
Since $f(a_1,\ldots,a_n)=f(a_{\sigma(1)},\ldots,a_{\sigma(n)})$
for any $\sigma \in \Sym(n)$, we are allowed to denote
$f(a_1,\ldots,a_n)$ by $f(a)$.

\section{Main Theorem}\label{S:main}

\begin{theorem} \label{main} Let $F$ be a pseudofinite field, $F(x)$ a field extension of $F$
with $x$ transcendental over $F$ and $H$ a nonabelian simple group.
Let $$g(T,Y_1,\ldots,Y_{n-1},Y_n)\in F[T,Y_1,\ldots,Y_{n-1},Y_n]$$
be a polynomial over $F$ symmetric in the indeterminates
$Y_1,\ldots,Y_n$. For every $a=\{a_1,\ldots,a_{n-1}\} \in
F^{[n-1]}$, let $L_{a}$ denote the splitting field of $g(T,a,x)\in
F(x)[T]$. Suppose that for every $a$, $b$ in $F^{[n-1]}$ the
following properties are satisfied:

\begin{enumerate}
    \item $\Gal(L_{a}/F(x))\cong H$.
    \item $L_{a}$ is a regular extension of $F$.
    \item $L_{a} \neq L_{b}$ for $a\neq b$.
\end{enumerate}

If $R \subset F^{[n]}$ is defined by the condition
$$\hbox{``}\{a_1,\ldots, a_n\}\in R \hbox{ if and only if   } \, g(T,a_1,\ldots,a_{n-1},a_n) \hbox{
has a root in } F \hbox{"} \leqno{(1)}$$ then $(F, R)$ is a random
$n$-hypergraph. \end{theorem}

We will construct polynomials satisfying the conditions stated in
the hypothesis after giving the proof of the theorem. This will
allow us to conclude that we can realize a random $n$-hypergraph in
a pseudofinite field $F$.

\vskip .1 in

\proof  Let $a^1=\{a^1_1,\ldots,a^1_{n-1}\},\ldots,
a^m=\{a^m_1,\ldots,a^m_{n-1}\}$ be in $F^{[n-1]}$ for $1 \leq
i\leq m$ and let $I \subseteq \{1,...,m\}$ and
$J=\{1,...,m\}\setminus I$. To prove that $R$ is a random
$n$-hypergraph we need to find $c\in F$ such that
$\{a^i_1,\ldots,a^i_{n-1},c\} \in R\,$ for every $i\in  I$, and
$\{a^j_1,\ldots,a^j_{n-1},c\} \not\in R \,$ for every $ j \in J$.

Unwinding the definition (1), this says that we need to find an
element $c$ of $F$ such that the polynomial
$g(T,a^i_1,\ldots,a^i_{n-1},c)$ will have a root in $F$ for  $i\in
I$  and the polynomial $g(T,b^j_1,\ldots,b^j_{n-1},c)$ will have
no roots in $F$ for $j\in  J$.

The strategy of the proof is as follows: we will construct an
elementary extension of $F$ containing an element $x$ satisfying the
conditions required for $c$. We can then conclude that such an
element exists in $F$ as well.

Let $L_1,\ldots ,L_m$ be the splitting fields of the polynomials
$$g(T,a^1_1,\ldots,a^1_{n-1},x),\ldots,
g(T,a^m_1,\ldots,a^m_{n-1},x)\in F(x)[T]$$ over $F(x)$ respectively.
These splitting fields are distinct extensions with non abelian
simple Galois groups by hypothesis, therefore they are linearly
disjoint by lemma \ref{setwise}. Denote by $L = L_1\ldots L_m$ the
join of the extensions $L_1,\ldots,L_m$. By Fact \ref{lindis},
$$\Gal(L/F(x)) \simeq H\times H\times \ldots \times H$$ is the
product of $m$ copies of $H$.

By Lemma \ref{transitive} there is an element $\mu_j$ of
$\Gal(L_j/F(x))$ which moves all the roots of
$g(T,a^i_1,\ldots,a^i_{n-1},x)$ for every $j\in J$. Take $\mu$ in
$\Gal(L/F(x))$ so that $\mu_{|L_i} = \Id$ for every $i\in I$ and
$\mu_{|L_j}=\mu_j$ for every $j\in J$.

Let $\sigma$ be a topological generator of the absolute Galois group
$\Gal(\widetilde{F}/F)\simeq \widehat {\mathbb{Z}}$ of the
pseudofinite field $F$. By the hypothesis of the theorem, the fields
$L_i$ are regular extensions of $F$. Since $L_1,\ldots,L_m$ are
linearly disjoint over $F(x)$, this implies that $L$ is a regular
extension of $F$ by Lemma \ref{reg}. Therefore there is an
automorphism $\tau\in \Gal(\widetilde{F(x)}/F(x))$ extending both
$\sigma$ and $\mu$. Denote the fixed field of $\tau$ by $M$.

Proposition \ref{M} implies that $M$ is a regular extension of $F$
and that it has a unique extension of degree $n$ for every $n\in
\mathbb{N}$. This condition with Theorem \ref{elementary} imply that
there is a pseudofinite field $E$ containing $M$ which is a regular
extension of $M$, i.e.\ $\tilde{M}\cap E =M$. Therefore $E$ is  a
regular extension of $F$ and so, by Theorem \ref{el}, $E$ is an
elementary extension of the pseudofinite field $F$ containing $x$.

Now we claim that
$$E\models \exists c ([\bigwedge_{i\in I} \exists T\,\, g(a^i_1,\ldots,a^i_{n-1},c,T)=0 ]\wedge
[\bigwedge_{j\in J} \forall T\,\,
g(a^j_1,\ldots,a^j_{n-1},c,T)\neq 0]).$$

Taking  $x\in F(x)< E$ for the variable $c$ in the above sentence,
we will prove, $$[\bigwedge_{i\in I}\,\exists T\,
g(a^i_1,\ldots,a^i_{n-1},x,T)=0 ]\wedge [\bigwedge_{j\in J}
\,\forall T\, g(a^j_1,\ldots,a^j_{n-1},x,T)\neq 0]$$ holds in the
pseudofinite field $E$.

Let $i\in I$,  $L_i$ contains all roots of
$g(a^i_1,\ldots,a^i_{n-1},x,T)$. $\mu$ is the identity on $L_i$ and
$\tau$ extends $\mu$, therefore $M$ the fixed field of  $\mu$,
contains all the roots of $g(a^i_1,\ldots,a^i_{n-1},x,T)$. The
pseudofinite field  $E$  contains $M$, so
$$E\models \bigwedge_{i\in I} \exists \,T\,\, g(a^i_1,\ldots,a^i_{n-1},x,T)=0.$$

We will show that $c=x$ satisfies the second part of the
conjunction. Suppose for a contradiction that for some $j \in J$
there exists a $t \in E$ such that
$g(a^j_1,\ldots,a^j_{n-1},x,t)=0$. Since $t$ is a root of the
polynomial $g(a^j_1,\ldots,a^j_{n-1},x,T)$, $t$ is in $L_j< L$, an
algebraic extension of $F(x)$, $t$ is also in $E$. Hence $t \in
E\cap \widetilde{F(x)} = M$ therefore $t\in M\cap L_j= \Fix(\mu)$.

But we chose $\mu$ so that it does not fix any root of
$g(a_j^1,\ldots,a_j^{n-1},x,T)$, a contradiction. We conclude that
$g(a^j_1,\ldots,a^j_{n-1},x,T)$ does not have any root in $E$ for
all $j\in J$. This proves our claim.

Hence, $$E\models \exists c ([\bigwedge_{i\in I} \exists T\,\,
g(a^i_1,\ldots,a^i_{n-1},c,T)=0 ]\wedge [\bigwedge_{j\in J}
\forall T\,\, g(a^j_1,\ldots,a^j_{n-1},c,T)\neq 0]).$$

The formula above has parameters from $F$. Since $E$ is an
elementary extension of $F$ it is also true that $$F\models \exists
c ([\bigwedge_{i\in I} \exists T\,\, g(a^i_1,\ldots,a^i_{n-1},c,T)=0
]\wedge [\bigwedge_{j\in J} \forall T\,\,
g(a^j_1,\ldots,a^j_{n-1},c,T)\neq 0]).$$ And this proves the
theorem. $\hfill{\Box}$

\section{Construction of Extensions}

In Section \ref{S:main} we proved Theorem \ref{main} which states
that, if there exists a polynomial $g(T,X_1,\ldots,X_n)$ over a
pseudofinite field $F$  satisfying certain conditions, then using
this polynomial we can define a random $n$-hypergraph on $F$. Here
in this section we will construct polynomials satisfying the
conditions of Theorem \ref{main} which will allow us to define a
random $n$-hypergraph on $F$.

The methods of constructing polynomials satisfying the conditions of
the Theorem \ref{main} vary with the characteristic of the given
pseudofinite field. We have two cases to consider separately:
characteristic 0, and positive characteristic. In both cases we will
use tools from the ramification theory of the function fields.

The following lemma describes the extensions of a valuation of a
function field $K(y)$ in an integral extension. It is an easy
consequence of \cite[Theorem III.3.7]{stich}.

\begin{lemma}\label{decomposition}\cite{stich, MM} Let $K(y)$ be a function field,
$f(X)\in K[X]$ a seperable monic polynomial and $g(X)=f(X)-y \in
K(y)[X]$. Let $\b \in K$ and $f(X)-\b=\prod \gamma_i(X)^{r_i}$ where
$\gamma_i$ are distinct irreducible polynomials in $K[X]$. Let $x$
be a root of $g(X)=f(X)-y$ and  $L=K(y)(x)$ be a root field of
$g(X)=f(X)-y$ over $K(y)$. Then the extensions of the valuation
$v_{y-\beta}$ of $K(y)$ to $L$ are the valuations $v_{\gamma_i(x)}$
and we have $r_i=r(v_{\gamma_i(x)}:v_{y-\beta})$.

\end{lemma}

The next lemma is an important result in the theory of valuations,
it can be found in \cite[Proposition III.8.9]{stich} stated in the
language of places instead of valuations.

\begin{fact}[Abhyankar's Lemma] \label{abhy} Let $L = L_1 L_2$
be the join of two finite algebraic extension fields $L_i$ of $K$.
Let $v$ be a valuation of $L$, whose restriction to $L_i$ is
ramified over $K$ with ramification index $r_i$. If at least one of
$r_i$ is not equal to $0$ modulo the characteristic of the field
$L$, then the ramification index of $v$ in $L/K$ is $\lcm(r_1,r_2)$.
\end{fact}

\subsection{Characteristic 0 Case}

We will work in characteristic $0$ throughout this section.

In \cite{serre} (p. 44) it is shown that the polynomial
$$h(T,Y)=(m-1)T^m-mT^{m-1}+1+(m-1)Y^2$$ gives rise to a regular Galois extension
of $\qq(Y)$ with Galois group equal to the alternating group on $m$
elements, $\Alt(m)$ when $m$ is divisible by 4. Let us denote the
splitting field of the polynomial $h(T,Y)$ over $\qq(Y)$ by $L$.
There are two valuations of $\qq(Y)$ which ramify in the extension
$L$: valuation $v_{\infty}$ of ramification index $m/2$ and the
valuation $v_{1+(m-1)Y^2}$ with ramification index $m-1$. Over
$\tilde{\qq}(Y)$ the valuation $v_{1+(m-1)Y^2}$ gives rise to two
valuations each with ramification index $m-1$.

We apply the linear transformation $Y \rightarrow 1/Y$ to the
polynomial $h(T,Y)$. This is a fractional linear transformation of
the base field hence $\Gal(h(T,Y),\qq(Y))=\Gal(h(T,1/Y),\qq(Y))$
Moreover we can multiply the resulting polynomial by $Y^2$ to
eliminate the denominators, this does not effect the Galois group
nor the ramification indices of the valuations. Therefore we obtain
the following lemma.

\begin{lemma}\label{polynomial} Let $m$ be a natural number divisible by 4.
Let $y$ be transcendental over $\qq$. The Galois group of the
polynomial $$f(T,y)=(m-1)y^2T^m-my^2T^{m-1}+y^2+(m-1)$$ over
$\mathbb{Q}(y)$ is the alternating group $\Alt(m)$. The polynomial
$f(T,y)$ is absolutely irreducible over $\qq$. There are two
valuations of $\mathbb{Q}(y)$ which ramify in the splitting field of
$f(T,y)$, the valuation $v_y$ with ramification index $m/2$ and the
valuation $v_{y^2+(m-1)}$ with ramification index $m-1$.
\end{lemma}


We fix a pseudofinite field $F$ of characteristic 0 and let $n\geq
3$ be a natural number and $x$ transcendental element over $F$. For
each $a\in F^{[n-1]}$ define $y_a=S(a_1,\ldots, a_{n-1},x)$, where
$S(t_1,\ldots,t_n)=\sum_{r_1+\ldots+r_n=n-1} t_1^{r_1}\ldots
t_n^{r_n}$ is the $(n-1)^{\hbox{st}}$ complete elementary symmetric
polynomial defined in section \ref{symmetric}. Then $y_{a}$ is
transcendental over $F$, $x$ is a root of the polynomial
$$S(a_1,\ldots, a_{n-1},X)-y_{a} \in F(y_a)[X]\leqno{(1)}$$
and $F(x)$ is a degree $n-1$ extension of $F(y_{a})$. By Fact
\ref{regular}, $F(x)$ is a regular extension of $F$. We call
$F(x)$ the ``small" extension of $F(y_{a})$.

Now we will build a Galois extension of $F(y_a)$ with Galois group
$\Alt(m)$.  Let $m=4(n-1)!$ and $K_{a}$ be the splitting field of
the polynomial
$$f(T,y_{a})=(m-1)y_{a}^2T^m-my_{a}^2T^{m-1}+y_{a}^2+(m-1)\leqno{(2)}$$
over $F(y_{a})$. By Lemma \ref{polynomial} the polynomial
$f(T,y_{a})\in\qq(y_a)[T]$ is absolutely irreducible over
$\mathbb{Q}$ and its Galois group over $\qq(y_{a})$ is $\Alt(m)$.
Hence by Fact \ref{stable}, the Galois group of $f(T,y_{a})$ over
$F(y_{a})$ is $\Alt(m)$  and $K_{a}$ is a regular extension of $F$.
We call $K_{a}$ the ``large" extension of $F(y_{a})$.

Note that if $m\geq 5$, $\Gal(K_{a}/F(y_{a}))=\Alt(m)$ is a simple
group. Since $m=4(n-1)!$, the degree of the extension
$[K_{a}:F(y_{a})]=m!/2$, is larger than $(n-1)!$ hence $K_{a}$
cannot be contained in the Galois closure of $F(x)$ over $F(y_{a})$
which is of degree at most $(n-1)!$. Then by Corollary \ref{simple},
we conclude that $K_{a}$ and $F(x)$ are linearly disjoint over $F$.

Now let $L_{a}$ be the join of $K_{a}$ and $F(x)$, the small and the
large extensions of $F(y_{a})$. Since $K_{a}$ and $F(x)$ are
linearly disjoint over $F(y_{a})$, $\Gal(L_{a}/F(x))\cong
\Gal(K_{a}/F(y_{a}))=\Alt(m)$. Also, since both $K_{a}$ and $F(x)$
are regular extensions of $F(y_{a})$, and since they are linearly
disjoint over $F(y_{a})$, their join $L_{a}$ is a regular extension
of $F$ by Lemma \ref{reg}.

Now consider the extension $L_a/F(x)$. Since $y_a=S(a,x) \in F(x)$
we see that $L_{a}$ is the splitting field of the polynomial
$$f(T,S(a,x))=(m-1)S(a,x)^2T^m-m
S(a,x)^2T^{m-1}+S(a,x)^2+(m-1)$$ over $F(x)$.

Note that $f(T,S(a,x))=g(T,a,x)$ where $g(T,X_1,\ldots X_n)$ is a
symmetric polynomial in $X_1,\ldots,X_n$ which is the desired
polynomial. Note that for every $a \in F^{n-1}$, the splitting field
$L_{a}$ of $g(T,a,x)$ over $F(x)$ (the field constructed above)
satisfies the following conditions of the Theorem \ref{main}: (i) $
\Gal(L_{a}/F(x))\cong \Alt(m)$, a simple non-abelian group since
$m=4(n-1)!>5$ for $n\geq 3$ (ii) $L_{a}$ is a regular extension of
$F$. Now we will prove condition (iii):

\vskip .1in

\claim (iii) For every $a$, $b$ in $F^{[n-1]}$ $L_{a} \neq L_{b}$ if
$a\neq b$.

\vskip .1in

\proof First note that $K_a$ and $L_{a}$ are regular extensions of
$F$ for every $a$ in $F^{[n-1]}$. Hence, working over  the algebraic
closure $\tilde{F}$ of $F$ instead of $F$ will not change the Galois
groups we have constructed.  We denote the extensions of
$\tilde{F}(x)$ that corresponds to the extensions $K_a$ and $L_a$ of
$F(x)$ by $\tilde{K}_a$ and $\tilde{L}_a$ for every $a$ in
$F^{[n-1]}$. To show that $L_a$ and $L_b$ are distinct extensions of
$F(x)$, it is enough to show that $\tilde{L}_a$ and $\tilde{L}_b$
are distinct extensions of $\tilde{F}(x)$.

We will find a valuation of the field $\tilde{F}(x)$ which has
different ramification indices in the Galois extensions
$\tilde{L}_a$ and $\tilde{L}_b$ hence conclude that $\tilde{L}_a\neq
\tilde{L}_b$ for $a\neq b$.

Let $a\neq b$ be in $F^{[n-1]}$,  we know by Lemma \ref{symmetric}
that $S(a,X)\neq S(b,X)$. Then there exists a factor $(X-\a)$ of the
polynomial $S(a,X) \in \tilde{F}[X]$ such that the multiplicity of
$(X-\a)$ in $S(a,X)$ is $e_1>0$ and the multiplicity of $(X-\a)$ in
the polynomial $S(b,X)$ is $e_2\geq 0$ where $e_1\neq e_2$.

Recall that by Lemma \ref{polynomial}  the valuations of $\qq(y)$
that ramify in the splitting field of
$f(T,y)=(m-1)y_{a}^2T^m-my_{a}^2T^{m-1}+y_{a}^2+(m-1)$ are $v_{y}$
and $v_{y^2+(m-1)}$ . Also recall that we denote the splitting field
of $f(T,y_{a})$ over $\tilde{F}(y_a)$ by $\tilde{K}_a$.  Since the
constant field $\tilde{F}$ of $\tilde{F}(y_a)$ is algebraically
closed, as a consequence of Lemma \ref{polynomial} the valuations of
$\tilde{F}(y_a)$  that ramify in the extension
$\tilde{K}_a/\tilde{F}(y_a)$  are $v_{y_a}$, $v_{y_a-\gamma}$ and
$v_{y_a+\gamma}$ where $\gamma$ is $(1-m)^{1/2}$. Moreover, since we
obtained the extension $\tilde{L}_a/\tilde{F}(x)$ by setting
$y_a=S(a,x)$ the valuations of $\tilde F(x)$  which ramify in
$\tilde L_a$ are exactly the valuations $v_{x-\beta}$, where $\beta$
is a root of $S(a,X)=0$, or of $S(a,X)=\pm\gamma$.

Now we will calculate the ramification index of the valuation
$v_{x-\a}$ of $\tilde{F}(x)$ in $\tilde{L}_a$ and $\tilde{L}_b$.

By Theorem \ref{decomposition}, the valuation $v_{x-\a}$ of
$\tilde{F}(x)$ extends the valuation $v_{y_a}$ of $\tilde{F}(y_a)$
with ramification index $r(v_{x-\a}:v_{y_a})=e_1$ since the
multiplicity of $X-\a$ in $S(a,X)$ is $e_1>0$. Also any extension of
the valuation $v_{y_a}$ to $K_a$ has ramification index $m/2$.

Since $e_1\leq n-1$, it divides $m/2=2(n-1)!$ and therefore
$r(w:v_{x-\alpha})=m/2e_1>1$ for any valuation $w$ on $\tilde L_a$
which extends $v_{x-\alpha}$.

Now we will calculate the ramification index of the valuation
$v_{x-\a}$ in the extension $L_b/F(x)$.

If $e_2>0$, the same argument shows that $v_{x-\alpha}$ ramifies in
$\tilde L_b$ with ramification index $m/2e_2\neq m/2e_1$.

If $e_2=0$, then either $v_{x-\alpha}$ does not ramify in $\tilde
L_b$, or else $S(b,\alpha)=\pm \gamma$, in which case, using
Abhyankar's Lemma, we obtain that the index of ramification of
$v_{x-\alpha}$ in $\tilde L_b$  divides $m-1$. In that case, since
$m$ and $m-1$ are relatively prime, and $m/2e_i>1$ divides $m$, we
also obtain that the indices of ramification of $v_{x-\alpha}$ in
$L_a$ and in $L_b$ are distinct. Since $\tilde{L}_a$ and
$\tilde{L}_b$ are Galois extensions of $\tilde{F}(x)$, this implies
that $\tilde{L}_a\neq \tilde{L}_b$. Hence $L_a\neq L_b$, and the
claim is proved.\qed

\vskip .2in

We showed that conditions (i),(ii),(iii) of  Theorem \ref{main} hold for the polynomial $g(T,X_1,\ldots,X_{n-1},X_n)$ for $n\geq 3.$ Hence, we
conclude that one can define a random $n$-hypergraph  in a pseudofinite field of characteristic 0 for $n\geq 3$. For $n=2$ the same method of
constructions can be applied by choosing $m=8$ to to satisfy the condition that $\Alt(m)$ is simple. This gives formula defining a random graph
in a pseudofinite field different from the one given by Duret.

\subsection{Positive Characteristic $p$: Enlarging the
Ramification Locus}

We will use  Abhyankar's polynomials to build Galois extensions of
function fields with Galois group $\Alt(m)$ in positive
characteristic. The following theorem gives us polynomials over
$K(y)$ with Galois group $\Alt(m)$ in case the characteristic of the
field is greater than 2.

\begin{theorem}[\cite{Ab0}, \cite{Ab1}]
\label{Alt3} Let $K$ be a field of characteristic $p>2$,  $y$
transcendental over $K$ and $L$ the splitting field of the
polynomial $f_{t,p}(T,y)=T^m-yT^t+1$ over $K(y)$ where $t\not \equiv
0 (\mod\, p)$ and $m=t+p$.  Then $\Gal(f_{p,t}(T,y),K(y)) \simeq
\Alt(m)$. Additionally the valuation $v_{\infty}$ of $K(y)$ splits
into the valuations $\hat{v}_{s}$ and $\hat{v}_{\infty}$ with
ramification indices $t$ and $p$ in the extension $K(s)$ of $K(y)$
where $s$ is a root of $f_{p,t}(T,y)=0$ over $K(y)$.
\end{theorem}

The following theorem is extracted from \cite{Ab2} where the above
result is extended for fields of characteristic two.

\begin{theorem}[\cite{Ab0}, \cite{Ab2}] \label{Alt2} Let $K$ be a field of
characteristic $p=2$,  $y$ transcendental over $K$ and $L$ the
splitting field of the polynomial $f_{t,q}(T,y)=T^m-yT^t+1$ over
$K(y)$ where $m=t+q$ and $t$ and $q$ satisfy the following
conditions:
\begin{enumerate}
    \item $q=p^l$ for some $l$
    \item $t\not \equiv 0 (\mod\, p)$
    \item $t>q>p$
    \item $t+q\equiv 1(\mod 8)$ or $t+q\equiv 7 (\mod 8)$

\end{enumerate}

Then $\Gal(f_{t,q}(T,y),K(y)) \simeq \Alt(m)$. Additionally the
valuation $v_{\infty}$ of $K(y)$ splits into the valuations
$\hat{v}_0$ and $\hat{v}_{\infty}$ with ramification indices $t$ and
$q$ in the extension $K(s)$ of $K(y)$ where $s$ is  a root of
$f_{p,t}(T,y)=0$ over $K(y)$.
\end{theorem}



Next, we apply the fractional linear transformation $y \mapsto 1/y$
of $K(y)$ to the polynomial $f_{t,q}(T,y)=T^m-yT^t+1$ to get the
polynomial $h_{t,q}(T,y)=yT^m-T^t+y$. Under this transformation the
valuation $v_{\infty}$ is sent to the valuation $v_0$. Combining
Theorems \ref{Alt3} and \ref{Alt2} with this transformation we have
the following corollary.

\begin{corollary}\label{Alt} Let $K$ be a field of
characteristic $p>0$,  $y$ transcendental over $K$ and $L$ the
splitting field of the polynomial $h_{t,q}(T,y)=yT^m-T^t+y$ over
$K(y)$ where $t\not \equiv 0 (\mod\, p)$ and $m=t+q$. Take  $q=p$ in
case the characteristic of the field $K$ is $p> 2$ and take $q=p^l$,
$t>q>p$  and $m\equiv 1(\mod 8)$ or $m\equiv 7 (\mod 8)$ if the
characteristic is $p=2$. Then $\Gal(h_{t,q}(T,y),K(y)) \simeq
\Alt(m)$. Additionally,  $v_y$ is the only valuation of $K(y)$ which
ramifies in $L$ and the ramification index of any extension of the
the valuation $v_{y}$ of $K(y)$ to $L$ is divisible by $t$.
\end{corollary}

Let $F$ be a pseudofinite field of positive characteristic $p$.
Let $n>1$ be such that $p\nmid n-1$. Let $x\in \Omega$ be
transcendental over $F$. We will construct polynomials satisfying
the conditions of Theorem\ref{main}.

For every $a=\{a_1,\ldots,a_{n-1}\}\in F^{[n-1]}$, let $z_a$ be
equal to $S(a_1\ldots,a_{n-1},x)$ where $x$ is the transcendental
element we fixed at the beginning. Then $F(x)$ is the field
extension of $F(z_a)$ given by the polynomial
$S(a_1\ldots,a_{n-1},X)-z_a$ of degree $n-1$ and $z_a\in \Omega$
is transcendental over $F$.

For $k>n-1$, let $p_1,\ldots,p_k$ be $k$ distinct primes greater
than $n-1$, not equal to the characteristic of $F$, each of which is
congruent to 1 or 7 modulo 8.  This condition is possible by
Dirichlet's theorem on arithmetic progression of primes. Also choose
$p_1,\ldots,p_k$ such that $p_1+\cdots +p_k$ is not congruent to 0
modulo the characteristic of the field.

Fix distinct $\b_1,\ldots,\b_k\in F$, we set $u_a=
\prod_{i=1}^k(z_a-\b_i)^{p_i}$. Then $F(z_a)$ is a separable
extension of $F(u_a)$ of degree $p_1+\ldots +p_k$  and $z_a$ is a
root of
$$\prod_{i=1}^k(Z-\b_i)^{p_i}-u_a\in F(u_a).$$

Now take $t$ to be $t=p_1\ldots p_k$. Since $p_i$'s were chosen to
be congruent to 1 or 7 modulo 8, \, the product $t=p_1\ldots p_k$ is
congruent to 1 or 7 modulo 8. Define $q$ by:
\begin{enumerate}
  \item $q=p$ if characteristic of $F$ is $p>2$,
  \item $q=8$ if the characteristic of $F$ is $2$.
\end{enumerate}

Let $m=q+t$ and $M_a$ be the splitting field of the polynomial
$f_{t,q}(T,u_a)=u_aT^m-T^t+u_a$ over $F(u_a)$. Note that, if
$\ch(F)=2$ then $m=8+t$ so $m$ is congruent to  1 or 7 modulo 8 and
if $\ch (F)=p\neq 2$ then $m=p+t$. Hence by Theorem \ref{Alt}, we
have $Gal(M_a/F(u_a))\simeq \Alt(m)$ for both cases, $\ch(F)=2$ and
$\ch(F)>2$. Also by Lemma \ref{acts}, $M_a$ is a regular extension
of $F$ since $f_{t,q}(T,u_a)$ is absolutely irreducible over $F$.

Note that $M_a$ is not contained in the Galois closure of $F(z_a)$
over $F(u_a)$ because $(t+q)!/2=(p_1\cdot\ldots\cdot
p_k+q)!/2>(p_1+\ldots+p_k)!$ and $\Gal(M_a/F(u_a))$ is a simple
group. Thus, by Corollary \ref{simple}, $M_a$ is linearly disjoint
from $F(z_a)$ over $F(u_a)$.

Let $K_a$ be the join of the extensions $M_a$ and $F(z_a)$. Then
$\Gal(K_a/F(z_a))=\Alt(m)$ since $M_a$  and $F(z_a)$ are linearly
disjoint over $F(u_a)$. Note that $K_a$ is the splitting field of
the polynomial
$$f_{t,q}(T,(\prod_{i=1}^k(z_a-\b_i)^{p_i}))=(\prod_{i=1}^k(z_a-\b_i)^{p_i})T^m-T^t+\prod_{i=1}^k(z_a-\b_i)^{p_i}\in F(z_a)[T]$$
over $F(z_a)$. And $K_a$ is a regular extension of $F$ by Fact
\ref{regular}.

So far we have constructed two extensions $K_a$ and $F(x)$ of
$F(z_a)$;  $K_a/F(z_a)$ Galois with Galois group $\Alt(m)$ and
$F(x)$ a finite algebraic extension of $F(z_a)$ given by the
polynomial $\prod_{i=1}^k(Z-\b_i)^{p_i}-u_a.$ Let $L_a$ be the join
of $F(x)$ with $M_a$. Note that $L_a$ is the splitting field of
$$f_{t,q}(T,\prod_{i=1}^k(S(a,x)-\b_i)^{p_i})=
\prod_{i=1}^k(S(a,x)-\b_i)^{p_i}T^m-T^t+\prod_{i=1}^k(S(a,x)-\b_i)^{p_i}$$
over $F(x)$.

Now let $g_1(T,a,x)=f_{t,q}(T, \prod_{i=1}^k(S(a,x)-\b_i)^{p_i}))$.
Then $g_1(T,X_1,\ldots,X_{n-1},X_n)$ is a symmetric polynomial in
$X_1,\ldots,X_n$. Note that (i) for every $a\in F^{[n-1]}$, the
Galois group of the polynomial $g_1(T,a,x)$ over $F(x)$ is the
simple group $\Alt(m)$, (ii) $L_a$ is a regular extension of $F$.
That is, the first two conditions for the main theorem are
satisfied. We also need to prove that

\,\,\, (iii) $L_a\neq L_b$ if $a\neq b$.

\proof The valuation $v_{u_a}$ of $F(u_a)$ extends to the valuations
$w_{z_a-\b_i}$ of $F(z_a)$ with ramification indices
$r(w_{z_a-\b_i}:v_{u_a})=p_i$ for $i=1\ldots,k$ by Lemma
\ref{decomposition}. Also the ramification index of the valuation
$v_{u_a}$ in the Galois extension $M_a$ is $s$ where $t|s$. Since
$p_i$ divides $t$, and $p_i$ is not divisible by the characteristic
of the field $F$, the valuations $w_{z_a-\b_i}$ of $F(z_a)$ have
extensions in $M_a$ which ramify with index $s/p_i$ over
$w_{z_a-\b_i}$ for $i=1,\ldots,k$ by Abhyankar's Lemma.

Again by Lemma \ref{decomposition}, the valuation $w_{z_a-\b_i}$ of
$F(z_a)$ extends to the valuations $v_{\gamma_1},\ldots,
v_{\gamma_h}$ of $F(x)$ according to the decomposition of the
polynomial $S(a_1\ldots,a_{n-1},X)-\b_i= \prod_{i=1}^h \gamma_i(X)$
in $F[X]$.

Since the characteristic of the pseudofinite field $F$ does not
divide $n-1$, the extension $F(x)/F(z_a)$ has at most $n-1$
ramification points. Therefore one of the valuations $w_{z_a-\b_i}$
does not ramify in $F(x)$ as $k$ was chosen to be greater than
$n-1$, fix one such $\b_i$ and let $\gamma_1(x),\ldots,\gamma_h(x)$
be the irreducible factors of $S(a,X)-\b_i$ over $F$. Then the
valuations $v_{\gamma_j(x)}$ for $1\leq j\leq h$ are the only
valuations of $F(x)$ extending $w_{z_a-\b_i}$ and they do not ramify
over $w_{z_a-\b_i}$ by our assumption on $\b_i$.

The ramification indices of the extensions $v_{\gamma_j(x)}$ of
$v_{z_a-\b_i}$ in $L_a$ are  $s/p_i$, by Abhyankar's Lemma. If
$j\neq i$, then the ramification indices over $F(x)$  of the
extensions of $v_{z_a-\beta_j}$ to $L_a$ divide $s/p_j$ and are
therefore different from $s/p_i$. It follows that we can retrieve
the polynomial $S(a,X)$ from the ramification locus of $L_a$ over
$F(x)$: choose $i$ such that the valuations of $F(x)$ ramifying with
index $s/p_i$ in $L_a$ are precisely
$v_{\gamma_1(x)},\ldots,v_{\gamma_r(x)}$, and $\prod_{j=1}^r
\gamma_j(X)$ has degree $n-1$ (such an $i$ exists by the discussion
above). Then $S(a,X)=\prod_{j=1}^k\gamma_j(X)-\beta_i$.

From this it follows that if $S(a,X)\neq S(b,X)$, then $L_a\neq
L_b$. This gives us the conclusion. \qed

We have showed the conditions of the main theorem are satisfied in a
pseudofinite field of characteristic $p$ where $p$ does not divide
$n-1$. Therefore we can interpret a random $n$-ary hypergraph in $F$
when the characteristic $p$ of $F$ is positive  and $p\nmid n-1$.
But if we can realize a random $m$-hypergraph, by restricting it to
$n<m$ many parameters, then we can realize random $n$-hypergraph as
well. Thus we can realize a random $n$-hypergraphs for every $n$.


\begin{thebibliography}{Abc1}
\bibitem[{\bf Ab}]{Ab0} S.\ Abhyankar, {\it Galois theory on the line in nonzero characteristic}, Bulletin of
American Mathematical Society \textbf{27} no 1 (1992) 68-133.

\bibitem[{\bf Ab1}]{Ab1} S.\ Abhyankar,
{\it Alternating group coverings of the affine line for
characteristic greater than two }, Ann.\ Math.\ \textbf{296}
(1993) 63-68.

\bibitem[{\bf Ab2}]{Ab2} S.\ Abhyankar, {\it Alternating group
coverings of the affine line for characteristic two},
Discrete Mathematics \textbf{133} (1994) 25-46.

\bibitem[{\bf Ax}]{Ax} J.\ Ax, {\it The elementary theory of finite fields},
Ann.\ Math.\ \textbf{88} (1968) 239-271.

\bibitem[{\bf B.et al}] {qm} L.\ B\'elair et al.,  {\bf Model theory and
Applications,} Quadranti di Mathematica, 2003.

\bibitem[{\bf Du}]{Duret} J.-L.\ Duret, {\it Les corps faiblement
alg\'ebriquement clos non s\'eparablement clos ont la propri\'et\'e
d'ind\'ependence}, in Model theory of Algebra and Arithmetic, Lect.\
Notes Math. \textbf{843} (1980), 135-157.

\bibitem[{\bf FJ}]{FJ} M.\ Fried, M.\ Jarden, {\bf Field Arithmetic}, Erg. Math. \textbf{11},
Berlin-Heidelberg-New York, 1986.

\bibitem[{\bf Fu}]{fulton} W.\ Fulton, {\bf Young Tableaux}, Oxford University Press, 1992.

\bibitem[{\bf Ho}]{Ho} W.\ Hodges, {\bf Model Theory}, Cambridge University Press, 1993.

\bibitem[{\bf H}] {pac}   E.\ Hrushovski, {\it Pseudo-finite fields and related structures}, Model Theory and Applications (eds. L. Bélair, et al.), 2003,
page

\bibitem[{\bf La1}]{lang} S.\ Lang, {\bf Algebra}, Adisson Wesley Publishing Company, 1993.
\bibitem[{\bf La2}]{lang2} S.\ Lang, {\bf Introduction to Algebraic Geometry}, Adisson Wesley Publishing Company, 1972.


\bibitem[{\bf MM}]{MM} G.\ Malle, B. H.\ Matzat, {\bf Inverse Galois Theory}, Springer-Verlag, Berlin, Heidelberg, New York, 1999.

\bibitem[{\bf Ma}]{marker} D.\ Marker, {\bf Model Theory: An Introduction},
Graduate texts in mathematics, Springer-Verlag, Berlin, Heidelberg,
New York, 2002.

\bibitem[{\bf Se}]{serre} J.P.\ Serre, {\bf Topics in Galois Theory}, Jones and Bartlett Publishers, 1992.

\bibitem[{\bf St}]{stich} H.\ Stichtenoth, {\bf Algebraic Function Fields and Codes},
Springer-Verlag , 1993.

\end{thebibliography}
\end{document}